\documentclass[12pt]{amsart}
\usepackage{amssymb}

\usepackage{helvet} % for Sans Serif fonts

\usepackage{courier} % default ttdefault

\reversemarginpar %marginal notes appear on the left side margin
\normalmarginpar %switches back to right side margin

\let\oldmarginpar\marginpar %changes the font of margin text
\renewcommand\marginpar[1]{\-\oldmarginpar{\raggedright\small\sf #1}}

\title{Restriction of sections of abelian schemes}

\author{Najmuddin Fakhruddin}
 
\address{School of Mathematics, Tata Institue of Fundamental Research, Homi Bhabha Road,
Mumbai 400005, India}

\email{naf@math.tifr.res.in}

\newcommand{\nc}{\newcommand}

\newcommand{\Q}{\mathbb{Q}}
\newcommand{\Z}{\mathbb{Z}}
\newcommand{\C}{\mathbb{C}}

\newcommand{\A}{\mathbb{A}}
\newcommand{\V}{\mathbb{V}}
\nc{\W}{\mathbb{W}}
\newcommand{\mc}{\mathcal}
\newcommand{\mb}{\mathbb}

\renewcommand{\P}{\mb{P}}

\nc{\aff}{{\A}^1}
\nc{\naive}{\!\sim_n}
\nc{\Spec}{\mathrm{Spec}}

\nc{\wt}{\widetilde}
\nc{\wh}{\widehat}
\newtheorem{thm}{Theorem}[section]
\newtheorem{prop}[thm]{Proposition}

\newtheorem{lem}[thm]{Lemma}

\theoremstyle{definition}

\newtheorem{rem}[thm]{Remark}

\input{xypic}
\begin{document}
\maketitle 

\section{Introduction}
\subsection{}

Let $A \to B$ be an abelian scheme. For a subvariety $C$ of $B$, we shall denote
by $A(C)$ the group of sections of the abelian scheme $A \times_B C \to C$.
We shall prove the following:
\begin{thm}
Let $B$ be a smooth, irreducible, quasi-projective variety over the complex
numbers and assume that $B$ has a  projective compactification $\bar{B}$
such that $\bar{B} - B$ is of codimension at least two in $\bar{B}$.
Then there exists a family of smooth ireducible curves $\{C_q\}_{q \in Q}$
in $B$ parametrised by an irreducible variety $Q$ such that
if  $p: A \to B$ is an abelian scheme
and  $q \in Q$ is a generic point,
then the restriction map on sections $A(B) \to A(C_q)$ is an isomorphism.
\label{thm:main}
\end{thm}

This answers, in a special case,  a question of Graber, Harris, Mazur and Starr
\cite[Question 4]{ghms2}.

\subsection{}
Our method of proof is briefly as follows: we first prove the
theorem for isotrivial abelian schemes and then reduce the general case to a cohomological
statement using the cycle class map. The Lefschetz hyperplane section
theorem allows us to further reduce to the case that $B$ is a surface.
This is then handled by a monodromy argument involving the
cohomology of Lefschetz pencils with coefficients in a local system.

\begin{rem}
It seems possible that our method can be extended to any smooth, quasi-projective
base $B$. However, the monodromy computations become much more difficult in this
generality. 
\end{rem}

\subsection{Conventions}
All our varieties will be over the field of complex numbers $\C$. By a generic
point of such a variety we shall mean a closed point lying outside a countable union
of proper closed subvarieties while by a general point we shall mean a closed point lying
in some Zariski open subset.

\section{Preliminary reductions}
\label{sec:reductions}
\subsection{} In this section $A \to B$, $A' \to B$ will always be  abelian schemes
with $B$ a smooth, irreducible quasi-projective variety of dimension $\geq 1$.
Unless stated otherwise, $C$ will be a smooth, irreducible curve in $B$ such that the map
$\pi_1(C) \to \pi_1(B)$ is surjective.
\begin{lem}
Suppose the map 
$A(B)\otimes \Q \to A(C) \otimes \Q$
is an isomorphism.
Then the map 
$A(B) \to A(C)$
is also an isomorphism.
\label{lem:tors}
\end{lem}

\begin{proof}
Since all elements of the  kernel of the restriction map are torsion, 
the kernel must be zero since a
non-zero torsion point always specializes to a non-zero torsion point.
For any $\sigma \in A(C)$ there exists $\tau \in A(B)$
and $n > 0$ such that
the restriction of $n \tau$ to $C$ is $\sigma$. 
Let $Z(\tau, n) = [n]^{-1}(\tau(B))$, where $[n]: A \to A$ is the multiplication
by $n$ map. $Z(\tau,n)$ is finite \'etale over $B$, so the surjectivity
assumption on fundamental groups implies that the restriction map
induces a bijection from the set of components of $Z(\tau, n)$ and those
of $Z(\tau,n)\times_B C$. Since $\sigma$ corresponds to a component
of $Z(\tau,n)\times_B C$ which is of degree $1$ over $C$, it follows that
it must be the restriction of an element of $A(B)$.
\end{proof}

\begin{lem}
Suppose the map
$A(B) \to A(C)$
is an isomorphism and let $A \to A'$ be an isogeny of abelian schemes over $B$.
Then the map 
$A'(B) \to A'(C)$
is also an isomorphism.
\label{lem:isog}
\end{lem}

\begin{proof}
This follows from
Lemma \ref{lem:tors} since $A(B) \otimes \Q \cong A'(B) \otimes \Q$.
\end{proof}

\begin{lem}
Suppose both the maps 
$A(B) \to A(C)$, and    $A'(B) \to A'(C)$
are isomorphisms.
Then the map 
$A \times_B A'(B) \to A\times_B A'(C)$
is also an isomorphism.
\label{lem:prod}
\end{lem}

\begin{proof}
This is clear since
$A \times_B A'(B) = A(B) \times A'(B) $.
\end{proof}

\begin{lem}
Suppose $A = A_0 \times B$ where $A_0$ is an abelian variety
i.e.~$A$ is a constant abelian scheme with fibre $A_0$. Let $\bar{B}$
be a normal projective compactification of $B$ and embed $\bar{B}$ in $\P^n$
for some $n$. Then for a generic
complete intersection curve $C$ in $B$ of large degree, the map
$A(B) \to A(C)$ is an isomorphism.
\label{lem:triv}
\end{lem}

\begin{proof}
Sections of a constant abelian scheme correspond to maps from the base to the fibre.
For a smooth, irreducible, quasi-projective variety $X$ we denote by $Alb(X)$
the Albanese variety of any smooth projective compactification of $X$. This is
universal for morphisms of $X$ to abelian varieties and is determined by
the mixed Hodge structure on $H_1(X,\Z)$ \cite{deligne-hodge3}.

Suppose $\dim(B) > 2$. Then by the  theorem 
of Goresky and MacPherson \cite[p. 150]{gm-morse}, it follows that if $B'$ is
a general hyperplane section of $B$ then the map $H_1(B',\Z) \to H_1(B, \Z)$
is an isomorphism. It follows that $Alb(B') \to Alb(B)$ is also an isomorphism
hence
\[
A(B) = Mor(Alb(B),A_0) \to Mor(Alb(B'),A_0) = A(B)
\]
is an isomorphism. 

We may thus assume that $\dim(B) = 2$. Let $\bar{C}$ be a general hypersurface
section of $\bar{B}$ and let $C = \bar{C} \cap B$. It follows from \emph{loc.~cit.}~that 
the map $H_1(\bar{C}, \Z) \to H_1(\wt{B}, \Z)$ is a surjection, where $\wt{B}$ is 
a resolution of singularities of $\bar{B}$.
Since $Alb(C) = Alb(\bar{C})$ and $Alb(B) = Alb(\wt{B})$, we get
an exact sequence of abelian varieties
\[
0 \to K \to Alb(C) \to Alb(B) \to 0 \ .
\]
If $\bar{B}$ is smooth, the usual theory of Lefschetz pencils \cite[Expos\'e XVIII]{sga72} 
implies that if $C$ is generic then $K$ is a
simple abelian variety which moreoever varies in moduli as $C$ varies, hence
$Hom(K,A_0) = 0$. If $B$ is only normal, using Lemma \ref{lem:lefs} one sees
that if the degree of hypersurface is sufficiently large then the usual arguments 
show that $K$ has no ``fixed part'' as $C$ varies so we still have
$Hom(K, A_0) = 0$ for generic $C$. Therefore $Hom(Alb(B), A_0) = Hom(Alb(C), A_0)$,
qnd so
\[
A(B) = Mor(Alb(B),A_0) \to Mor(Alb(C),A_0) = A(C)
\]
is an isomorphism.
\end{proof}

\begin{lem}
Suppose $A \to B$ is isotrivial
(i.e.~there exists a finite \'etale cover $B' \to B$ such that
the abelian scheme $A\times_B B' \to B'$ is constant) and let 
$C$ be a generic complete intersection
curve of large degree. Then the map
$A(B) \to A( C)$
is an isomorphism.
\label{lem:isotriv}
\end{lem}

\begin{proof}
We first note that an isotrivial abelian scheme is determined by
its monodromy representation $\pi_1(B) \to Aut(A_0)$, where
$A_0$ is the fibre over the basepoint.

Suppose $A(B)$ is non-torsion. 
Let $A''$ be the connected component containing the zero section 
of the Zariski closure of the union of the images of the 
elements of $A(B)$.
This is a (non-trivial) abelian subscheme of $A$ and is moreover
constant since the monodromy acts trivially on it.
Then $A$ is isogenous to $A'' \times_B (A/A'')$, with $A/A''$ also isotrivial.
The lemma follows in this case 
from Lemmas \ref{lem:triv} and \ref{lem:prod} along with induction
on the relative dimension.

If $A(B)$ is torsion, then 
$A( C)$ must also be torsion. For otherwise
$A \times_B C$ would contain a non-trivial constant abelian subscheme (by the argument above)
which is not possible because of the surjectivity of $\pi_1(C) \to \pi_1(B)$.
So the lemma follows from
Lemma \ref{lem:tors}. 
\end{proof}

\subsection{}
\label{sec:cohred}
Let $p:A \to B$ be an abelian scheme of relative dimension $n$ and consider the local system
$\V = R^{2n-1}p_*\Q$  on $B$. The monodromy representation
corresponding to this local system is semisimple, so breaks up as a direct
sum of irreducible representations. Let $\V = \V_1 \oplus \V_2$, where
$\V_1$ is the local system 
corresponding to the direct sum of all the irreducible summands with finite
image and $\V_2$ corresponds to the direct sum of those with infinite image.
Then there exist abelian schemes $p_i:A_i \to B$, $i = 1 ,2$
such that $A$ is isogenous to $A_1 \times_B A_2$ and 
$R^{2n-1}{p_i}_*\Q \cong \V_i$. $A_1$ is isotrivial and $A_2$ contains no isotrivial
abelian subschemes. By
Lemmas \ref{lem:isog}, \ref{lem:prod} and \ref{lem:isotriv}
the  proof of Theorem \ref{thm:main} reduces to the 
case that $p:A \to B$ is an abelian scheme which has no isotrivial abelian
subschemes.

Now suppose that $p:A \to B$ is as above with $ \V_1 = 0$.
Recall from \cite{saito-mhm} (or see \cite{arapura})
that the vector space $H^1(B, \V)$
carries a natural mixed Hodge structure of weights $\geq 2n$.
\begin{lem}
 $A(B)$ is a finitely generated abelian group 
and the cycle
class map identifies\footnote{For the proof of 
Theorem \ref{thm:main} the surjectivity of the map is not essential.} $A(B)\otimes \Q$ 
(via the Leray spectral sequence)
with the group of Hodge classes of
type $(n,n)$ of $H^1(B, \V)$. 
\label{lem:leray}
\end{lem}

\begin{proof}
The finite generation follows from the theorem of Lang and N\'eron  \cite{lang-neron}.
\smallskip

To prove the injectivity of the map we may assume that $B$ is a curve.
The cycle class map referred to above (after tensoring with $\Q$) 
can be viewed as being induced by the boundary map
\[
H^0(B, A) = A( B)_{an} \stackrel{\partial}{\to} H^1(B, \V_{\Z})
\]
coming from the long exact cohomology sequence corresponding to the short exact sequence
of sheaves on $B$ in the \emph{analytic} topology
\[
0 \to \V_{\Z} \to Lie(A/B) \to A \to 0 \ .
\]
Here
 $A(B)_{an}$ denotes the abelian group
of complex analytic sections,
$\V_Z = R^{2n-1}p_*\Z$ and $Lie(A/B)$, the relative Lie algebra of $A$ over $B$,
is the locally free sheaf associated to the complex local system $\V_{\Z} \otimes \C$.

Let $\bar{B}$ be the smooth compactification of $B$ 
and let $\bar{p}: \bar{A} \to \bar{B}$
be the N\'eron model of $A$. If $\partial(\sigma) = 0$ for some
$0 \neq \sigma \in A(B)$, then it lifts to a
non-zero element $\wt{\sigma}$ of $H^0(B, Lie(B/A))$. Since this is a complex
 vector space, 
the images of the elements $t \wt{\sigma}$, $t \in \C^*$ give 
a $1$-parameter family of elements of $A( B)_{an}$.
 By the property of Neron models, $\sigma$ extends to an element 
$\bar{\sigma} \in \bar{A}(\bar{B})$, hence by continuity any element
of $A( B)_{an}$ close to $\sigma$ also extends to an element of
$\bar{A}(\bar{B})_{an} = \bar{A}(\bar{B})$.
This implies that $A$ must contain a non-trivial 
constant abelian subscheme, which contradicts the hypotheses.
\smallskip

The surjectivity follows from Lefschetz's theorem on $(1,1)$ classes
since cup product with $\theta^{n-1}$, where $\theta  \in H^0(B, R^2p_*\Q)$
is the class of a polarisation,
induces an isomorphism of mixed Hodge structures
\[
H^1(B, R^1p_*\Q) \to H^1(B, R^{2n-1}p_*\Q) \otimes \Q(n-1) \ .
\]
\end{proof}

\begin{rem}
If $C$ is a generic complete intersection
curve in $B$ then the restriction map
$H^1(B, \V) \to H^1(C, \V|_C) $ is always
an injection; the difficulty lies in showing that all Hodge classes
of type $(n,n)$ lie in the image.
\end{rem}

\section{Lefschetz pencils with coefficients}

\label{sec:lefs}
\subsection{}
In this section we  state  mild generalisations of a couple of the results of the theory
of Lefschetz pencils.

\begin{lem} 
Let $D$ be the open unit disc in $\C$, $X$ a connected, two dimensional,
complex manifold and $\pi:X \to D$
a proper analytic map whose differential is non-zero except at a single point $x_0 \in X$
above $0 \in D$ where it has a non-degenerate critical point. Let $\V$ be a local
system of $\Q$-vector spaces on $X$. Then the monodromy of $R^1\pi_*\V$ restricted
to $D^* = D - \{0\}$ is unipotent.
\label{lem:unip}
\end{lem}

\begin{proof}
This follows from the results in SGA7 II \cite[Expos\'es XIII \& XIV]{sga72}:
we only indicate the slight changes that need to be made.
 Let $y \in D^*$ be a basepoint and 
$\sigma: H^1(X_y,\V) \to H^1(X_y,\V)$ be the monodromy automorphism.
Then $\sigma - 1$ is the composite 
of the following sequence of maps 
\[
H^1(X_y, \V) \stackrel{j^*}{\to}  H^1(V, \V) \stackrel{\mathrm{Var}}{\longrightarrow} 
H^1_c(V^0, \V) \stackrel{j^0*}{\to} H^1(X_y, \V)
\]
where $V$ is defined on \cite[p.~136]{sga72},
$j:V \to X_y$, $j^0:V^0 \to X_y$,  are  the inclusions and, following \cite[p.~151]{sga72},
  the map $\mathrm{Var}$ is as follows:
since $\V$ is constant on $V$ we may write $\V|_V \cong \oplus_{i=1}^r\Q\cdot v_i$,
where $v_i$, $i = 1,2, \ldots, r=rank(\V)$, is a local basis of sections. This gives rise
to elements $\delta_i \in H^1_c(V^0, \V)$ which are well defined upto sign. The basis
also gives an isomorphism $\V \to \check{\V}$ restricted to $V$, where $ \check{\V}$
is the dual local system of $\V$, and hence induces an isomorphism $H^1_c(V^0, \V) \cong
H^1_c(V^0,\check{\V} )$. We let $\check{\delta}_i \in H^1_c(V^0,\check{\V})$,
$ = 1,2, \ldots, r$ be the elements corresponding to the $\delta_i$'s under this isomorphism.
Then for any $x \in  H^1(V, \V)$
\[
\mathrm{Var}(x) = - \sum_{i=1}^r (x \cdot\check{\delta}_i )\, \delta_i \ ,
\]
where the pairing $(x \cdot \check{\delta}_i)$ is the natural duality pairing.

It is clear that $(j^*(j^0_*(\delta_k))\cdot \check{\delta}_l) = 0$ for all
$k,l = 1,2, \ldots, r$, hence $(\sigma -1)^2 = 1$. Thus $\sigma$ is unipotent.
\end{proof}

\begin{lem}
Let $X$ be a normal projective surface, 
$S$ a finite subset of $X$ including all its singular points
and $Y = X - S$.
There exists a pencil of very ample curves on $X$, $\{C_p\}_{p \in \P^1}$,
with the following properties:
\begin{enumerate}
\item A general element of the pencil is smooth and 
all $C_p$, $p \neq \infty$ are irreducible with at most a single
ordinary double point.
\item $S \subset C_{\infty}$
\item $C_p$ and $C_q$ for $p \neq q$ meet transversally (at smooth points of $X$).\\
\end{enumerate}
\label{lem:lefs}
\end{lem}

\begin{proof}
Let $\mc{L}$ be a very ample line bundle on $X$ and choose a trivialisation
of $\mc{L}$ restricted to $S$ \emph{i.e.}~an isomorphism $\mc{L}|_S \cong 
\mc{O}_S$ (which induces a similar isomorphism for all tensor powers of $\mc{L}$).
For $n > 0$, let $V_n$ be the subspace of $H^0(X, \mc{L}^{\otimes n})$
consisting of all sections whose restriction to $S$ is a constant section.
If $n$ is sufficiently large, this linear system is base point free and induces a
morphism $\phi_n:X \to \P(V_n)$ which is an embedding on $X - S$ and maps
$S$ to a single point.

Fix $n$ as above and let $X_n = \phi_n(X) \subset \P(V_n)$.
Let $\check{X_n} \subset  \check{\P}(V_n)$ be the dual variety of $X_n$.
It consists of two irreducible components, the general point of one 
corresponding to hyperplanes in $\P(V_n)$ tangent to a smooth point
of $X_n$ and the points of the other corresponding to hyperplanes containing
$\phi_n(S)$.

The proof of the existence of Lefschetz pencils in \cite{sga72} Expos\'e XVII
goes through without any changes
to show that a general pencil in $\check{\P}(V_n)$ gives rise to 
a pencil of hyperplane sections of $X_n$, which when pulled back to $X$ satisfies all
the conditions of the lemma. 
\end{proof}

\section{Proof of the theorem}
\label{sec:main}
\subsection{}
Let $p:A \to B$ be an abelian scheme of relative dimension $n$ with $B$ a smooth
connected surface and let $\V = R^{2n-1}p_*\Q$. Let $Y = B$
and let $X$ be a normal projective compactification of $Y$ with $S = X - Y$ 
a finite set. Assume that $\V$ 
does not contain any non-trivial sub-local systems with finite monodromy.
\begin{prop}
Let $\mc{L}$ be an ample line bundle on $X$.
Then for $n$ sufficiently large (depending
only on $X$) and $C$ a generic element of $\big |H^0(X, \mc{L}^{\otimes n})\big |$,
the retriction map
$A(B) \to A( C)$
is an isomorphism.
\label{prop:surf}
\end{prop}

\begin{proof}
We have seen in Section \ref{sec:cohred} that it suffices to prove that the 
restriction map $H^1(Y, \V) \to H^1(C, \V|_C)$
induces a surjection on Hodge classes of type $(n,n)$. 
\smallskip

Let $f: \wt{X} \to \P^1$ be obtained by blowing up the basepoints of a
pencil on $X$ obtained by applying Lemma \ref{lem:lefs}. Let $\wt{Y}$ be
the inverse image of $Y$ in $\wt{X}$ and let $\wt{\V}$ be the pullback of $\V$ to $\wt{Y}$.
Let $U$ be the open subset of $\P^1$ over which $f$ is smooth,
let $Y' = f^{-1}(U)$,  $f' = f|_{Y'}: Y' \to U$, $\V' = \wt{\V}|_{Y'}$
and let $\W = R^1{f'}_* \V'$. 
Our hypothesis on the monodromy of $\V$
implies that ${f'}_* \V' = 0$, hence from the Leray spectral sequence for $f'$
it follows that $H^0(U, \W) = H^1(Y', \V')$.

Let $u$ be a generic element of $U$ and let $C = C_u = f^{-1}(u)$.
Let $\alpha$ be a Hodge class of type $(n,n)$ in $H^1(C, \V|_C)$.
  Since $u$ is
generic and $\alpha$ is an algebraic class, it follows that there 
exist an \'etale morphism $g:T \to U$ with $T$ a smooth connected curve,
$t \in T$ such that $g(t) = u$,
and an element $\beta \in H^0(T, g^* \W)$
 such that $t^*(\beta) = \alpha$. 
It follows from Lemma \ref{lem:unip} that the local monodromies
of $\W$ around all the points of $\P^1 - (U \cup \{\infty\})$ are unipotent.
Since the image of $\pi_1(T) $ in $\pi_1(U)$ is of finite index and
$\A^1 = \P^1 - \{\infty\}$ is simply connected, it follows
that $H^0(U, \W) = H^0(T, g^*\W)$. So $\beta$ is already
defined over $U$, hence comes from an element of $H^1(Y', \V')$.

The proof is completed by observing that the pullback map $H^1(Y,\V) \to 
H^1(\wt{Y}, \wt{\V})$ is an isomorphism and that the restriction map
$H^1(\wt{Y}, \wt{\V}) \to H^1(Y', \V')$ induces a surjection 
(in fact an isomorphism) on 
Hodge classes of type $(n,n)$\footnote{This is elementary in our situation
since we know that such classes are algebraic.}.
\end{proof}

\subsection{}
We now complete the proof of the main result of the paper:
\begin{proof}[Proof of Theorem \ref{thm:main}]
Choose a normal, projective compactification $\bar{B}$ of $B$ with
$\bar{B} - B$ of codimension at least two in $\bar{B}$ and embed $\bar{B}$
in $\P^n$ for some $n$.
If $B'$ is a general hypersurface 
section of $B$ a
theorem of Goresky and MacPherson \cite[p.~150]{gm-morse}
implies that if $\dim(B') > 1$ then
$\pi_1(B') \to \pi_1(B)$ is an isomorphism,
and so also the restriction map 
\[
H^1(B, R^{2n-1}p_*\Q) \to H^1(B',R^{2n-1}p_*\Q) \ .
\]
 Since a complete intersection
curve in $B'$ is also a complete intersection curve in $B$, it follows from
the above and the reductions in Section \ref{sec:reductions}, that
 it suffices to prove the theorem  in the case $B$ is a surface.
This follows  from the reductions in Section \ref{sec:reductions} 
and Proposition \ref{prop:surf}. We see that the family of curves $\{C_q\}_{q \in Q}$
can be chosen to be all smooth complete intersection curves in $B$ of
a fixed  large multidegree.
\end{proof}

\subsection{Acknowledgements}
We thank Chandrashekhar Khare for a lively correspondence which provided the impetus
for writing this note and V.~Srinivas for a useful suggestion concerning the proof of 
Lemma \ref{lem:lefs}. We are grateful to Madhav Nori for pointing out some errors in
the first draft.

\end{document}